\documentclass[12pt]{amsart}
\usepackage{amsmath,amsthm,amsfonts,amssymb,eucal, hyperref}

\newcommand{\<}{\langle}
\renewcommand{\>}{\rangle}

\newcommand{\oq}{\ {\raise 7pt\hbox{${\scriptstyle\circ}$}}
\kern -7pt{
\hbox{$Q$}}}

\newcommand{\R}{ \mathbb R}

\newcommand{\C}{ \mathbb C}
\newcommand{\N}{ \mathbb N}

\newcommand {\ba}{\mathbf a}

\newcommand {\be}{\mathbf e}

\newcommand {\bk}{\mathbf k}

\newcommand {\bm}{\mathbf m}

\newcommand {\bze}{\mathbf 0}

\newcommand {\bn}{\mathbf n}

\DeclareMathOperator{\range}{{range}}

\DeclareMathOperator{\mb}{{mb}}
\DeclareMathOperator{\diag}{{diag}}

\newcommand{\1}
{{\,\vrule depth3pt height9pt}{\vrule depth3pt height9pt}
{\vrule depth3pt height9pt}{\vrule depth3pt height9pt}\,}

\DeclareMathOperator {\dist} {{dist}}

\DeclareMathOperator {\diam} {{diam}}

\DeclareMathOperator{\reg}{{reg}}
\DeclareMathOperator{\sing}{{sing}}

\DeclareMathOperator{\ran}{{ran}}

\newtheorem{thm}{Theorem}[section]
\newtheorem{cor}[thm]{Corollary}
\newtheorem{cla}[thm]{Claim}
\newtheorem{lem}[thm]{Lemma}
\newtheorem{prop}[thm]{Proposition}

\theoremstyle{definition}
\newtheorem{defn}[thm]{Definition}

\newtheorem{rem}[thm]{Remark}

\numberwithin{equation}{section}

\newcommand{\bee}{\begin{equation}}
\newcommand{\ene}{\end{equation}}
\newcommand{\bees}{\begin{equation*}}
\newcommand{\enes}{\end{equation*}}
\newcommand{\bes}{\begin{split}}
\newcommand{\ens}{\end{split}}

\newcommand{\bet}{\begin{thm}}
\newcommand{\ent}{\end{thm}}
\newcommand{\bel}{\begin{lem}}
\newcommand{\enl}{\end{lem}}
\newcommand{\bec}{\begin{cor}}
\newcommand{\enc}{\end{cor}}
\newcommand{\becl}{\begin{cla}}
\newcommand{\encl}{\end{cla}}
\newcommand{\bep}{\begin{proof}}
\newcommand{\enp}{\end{proof}}
\newcommand{\ber}{\begin{rem}}
\newcommand{\enr}{\end{rem}}
\newcommand{\ep}{\varepsilon}

\newcommand{\Z}{\mathbb Z}

\makeatletter
\def\square{\RIfM@\bgroup\else$\bgroup\aftergroup$\fi
  \vcenter{\hrule\hbox{\vrule\@height.6em\kern.6em\vrule}\hrule}\egroup}
\makeatother

 \usepackage{hyperref}

\def\dio{\mathrm{dio}}

\def\range{\mathtt{Range}}

\def\be{\mathbf{e}}
\newcommand{\one}{\mathbf{1}}

\usepackage[margin=1.1in]{geometry}

\begin{document}

\title[Perturbation theory]
{Perturbation theory}
\title[Convergence of perturbation theory]
{Perturbative diagonalisation for Maryland-type quasiperiodic operators with flat pieces}
\author[I. Kachkovskiy]
{Ilya Kachkovskiy}
\address{Department of Mathematics\\ Michigan State University\\
Wells Hall, 619 Red Cedar Rd\\ East Lansing, MI\\ 48824\\ USA}
\email{ikachkov@msu.edu}
\author[S. Krymskii]
{Stanislav Krymski}
\address{Department of Mathematics and Computer Science, St. Petersburg State University, 14th Line 29B, Vasilyevsky Island, St. Petersburg 199178, Russia}
\email{krymskiy.stas@yandex.ru}
\author[L. Parnovski]
{Leonid Parnovski}
\address{Department of Mathematics\\ University College London\\
Gower Street\\ London\\ WC1E 6BT\\ UK}
\email{leonid@math.ucl.ac.uk}
\author[R. Shterenberg]
{Roman Shterenberg}
\address{Department of Mathematics\\ University of Alabama, Birminghan\\
Campbell Hall\\1300 University Blvd\\ Birmingham, AL\\ 35294\\USA }
\email{shterenb@math.uab.edu}



\date{\today}



\begin{abstract}
We consider quasiperiodic operators on $\Z^d$ with unbounded monotone sampling functions (``Maryland-type''), which are not required to be strictly monotone and are allowed to have flat segments. Under several geometric conditions on the frequencies, lengths of the segments, and their positions, we show that these operators enjoy Anderson localization at large disorder.
\end{abstract}
\maketitle
\section{Introduction}
This paper can be considered as a direct continuation of the earlier publication \cite{KPS}. We consider quasiperiodic Schr\"odinger operators on $\ell^2(\Z^d)$:
\bee
\label{eq_h_def}
(H(x)\psi)_\bn=\varepsilon(\Delta\psi)_{\bn}+f(x+\omega\cdot\bn)\psi_{\bn},
\ene
where $\omega=(\omega_1,\ldots,\omega_d)$ is the {\it frequency vector}, and $\Delta$ is the discrete Laplacian:
$$
(\Delta\psi)_{\bn}=\sum_{\bm\in \Z^d\colon |\bm-\bn|_1=1}\psi_\bm.
$$
We will consider the regime of large disorder which, after rescaling, corresponds to small $\ep>0$.
The function $f$, which generates the quasiperiodic potential, is a non-decreasing function
\bee
\label{eq_f_infinity}
f\colon(-1/2,1/2)\to (-\infty,+\infty),\quad f(-1/2\pm 0)=\mp \infty,
\ene
and is extended into $\R\setminus(\Z+1/2)$ by $1$-periodicity. As usual for quasiperiodic operators, the numbers $\{1,\omega_1,\ldots,\omega_d\}$ are assumed to be linearly independent over $\mathbb Q$. 

Potentials of the form \eqref{eq_f_infinity} will be called {\it Maryland-type}, after the classical Maryland model with $f(x)=\tan(\pi x)$ (we refer the reader to the following, certainly not exhaustive, list of related publications \cite{Bellissard,FP,GFP,Simon_maryland,Wencai,JK,JY2,JY,Ilya}). In \cite[Section 6]{KPS}, we considered a case where $f$ was not strictly monotone and was allowed to have a single flat piece of sufficiently small length $h$; in particular, $|h|<\min_j |\omega_j|$. In this situation, we can treat this flat piece as a single isolated resonance. In the present paper, we consider several more elaborate situations where the resonance is not isolated but has, in some sense, finite multiplicity.

The main result of the paper is a general Theorem \ref{th_main_abstract}, where we establish some sufficient conditions under which the operator $H$ admits Anderson localization. We also list several particular examples and refinements in Section 6. The simplest new example (see Theorem \ref{th_main_1d} for a precise statement) is as follows. On $\ell^2(\Z)$, consider the operator \eqref{eq_h_def} with the function $f$ satisfying
$$
f(x)\equiv E,\quad x\in [a,a+L]\subset (-1/2,1/2).
$$
Additionally, assume that $\omega$ is Diophantine (see Section 2.2 for the definition), and that $[a-2\omega,a+L+2\omega]\subset (-1/2,1/2)$. Suppose that $f$ is Lipschitz monotone outside $[a,a+L]$, with some regularity conditions similar to \cite{KPS}. Assume also that $L$ is not a rational multiple of $\omega$. Then the operator \eqref{eq_h_def} has Anderson localization for $0<\ep<\ep_0(\omega,L)$.


We would also like to mention the simplest possible class of operators that are not covered by our approach. Suppose that $f$ is constant on an interval $[a,a+L]\subset (-1/2,1/2)$, and suppose that the set
$$
S:=\{\bn\in\Z^d\colon x+\omega\cdot\bn\in [a,a+L]+\Z\}
$$
has an unbounded connected component for some $x$ (here we define connectedness on the $\Z^d$ graph with nearest neighbour edges). For example, this will happen if $0<\omega_1<\omega_2<L$. Our methods cannot cover such operators. In fact, it seems possible that such models do not demonstrate Anderson localization around the energy $E=\left.f\right|_{[a,a+L]}$. 

Our approach, essentially, is based on assuming the opposite of the above: that is, that all connected components of the set $S$ are bounded and are sufficiently far away from each other. If we are able to surround each component by a layer of lattice points $\bn$ such that $f(x+\omega\cdot\bn)$ has good monotonicity properties, then, after a partial diagonalisation, their monotonicity  will ``propagate'', in a weaker form, into the interior of the set $S$. Afterwards, we can apply a modified version of the main result of \cite{KPS} to finalize the diagonalisation, see Proposition \ref{prop_convergence}.
\subsection*{Acknowledgments} The authors would like to dedicate this paper to the memory of Jean Bourgain.

The research of LP was partially supported by EPSRC grants EP/J016829/1 and EP/P024793/1. RS was partially supported by NSF grant DMS--1814664. IK was partially supported by NSF grant DMS--1846114.
\section{Preliminaries: regularity of $f$ and convergence of perturbation series}
While the functions $f$ under consideration will have flat pieces, a necessary assumption for all proofs below would be the existence of sufficiently many pieces with a good control of monotonicity. The corresponding regularity conditions and convergence results are summarized in this section.

It will be convenient to {\it not} exclude the case $x\in 1/2+\Z+\omega\cdot\Z^d$, where exactly one value of the potential in the operator $\eqref{eq_h_def}$ becomes infinite: say, $f(x+\omega\cdot\bn)$. In this case, the natural limiting object is the operator on $\ell^2(\Z^d\setminus\{\bn\})$ obtained from \eqref{eq_h_def} by enforcing the Dirichlet condition $\psi(\bn)=0$. The results of \cite{KPS}, which we will be using, extend ``continuously'' into these values of $x$, see \cite[Remark 4.14]{KPS}, with a reasonable interpretation of infinities, if one adds an infinite eigenvalue with an eigenvector $e_{\bn}$. Here, $\{e_{\bn}\colon \bn\in \Z^d\}$ is the standard basis in $\ell^2(\Z^d)$. We will also use the notation $\{\be_j:j=1,2,\ldots,d\}$ for the standard basis in $\R^d$. For a subset $A\subset \R^d$, it will be convenient to use the notation $\ell^2(A)$ for $\ell^2(A\cap \Z^d)$. For example, $\ell^2([0,N])=\ell^2\{0,1,2,\ldots,N\}$.

\subsection{$C_{\reg}$-regularity} Similarly to \cite{KPS}, we will always assume the following:
\begin{itemize}
	\item[(f1)] $f\colon (-1/2,1/2)\to \R$ is continuous, non-decreasing, $\quad f(-1/2+0)=-\infty,\quad f(1/2-0)=+\infty$, and $f$ is extended by $1$-periodicity into $\mathbb R\setminus(\Z+1/2)$.
\end{itemize}
Suppose, $f$ satisfies $\mathrm{(f1)}$. Let $C_{\mathrm{reg}}>0$, and $x_0\in (-1/2,1/2)$. We say that $f$ is $C_{\mathrm{reg}}$-regular at $x_0$, if:
\begin{itemize}
\item[(cr0)]The pre-image $f^{-1}((f(x_0)-2,f(x_0)+2))\cap (-1/2,1/2)$ is an open interval (denoted by $(a,b)$), and $\left.f\right|_{(a,b)}$ is a one-to-one map between $(a,b)$ and $(f(x_0)-2,f(x_0)+2)$.
\item[(cr1)]Let $D_{\min}(x_0):=\inf\limits_{x\in (a,b)}f'(x)\ge 1$ (for points where $f'$ does not exist, consider the smallest of the derivative numbers). Then,
\bee
\label{eq_cr1}
D_{\min}(x_0)\le f'(x)\le C_{\mathrm{reg}} D_{\min}(x_0),\quad \forall x\in (a,b).
\ene
\item[(cr2)]Define $(a_1,b_1)=f^{-1}(f(x_0)-1,f(x_0)+1)\subset (a,b)$, and

$$
g(x)=\frac{1}{f(x)-f(x_0)}, \quad x\in (b_1,a_1+1),
$$
extended by continuity to $g(\pm 1/2)=0$ (recall that we also assume $f(x+1)=f(x)$, so that the interval $(b_1,a_1+1)$ is essentially $(-1/2,1/2)
\setminus(a_1,b_1)$ together with the point $1/2=-1/2\,\,\mathrm{mod}\,\,1$). Then, under the same conventions on the existence of derivatives,
$$
|g'(x)|\le C_{\mathrm{reg}} D_{\min}(x_0),\quad x\in (b_1,a_1+1).
$$
\end{itemize}
For convenience, we will require $D_{\min}(x_0)\ge 1$ as a necessary condition for $C_{\reg}$-regularity. This condition can always be achieved by rescaling.
\subsection{The frequency vector}
The frequency vector $\omega\in \R^d$ is called {\it Diophantine} if there exist $C_{\dio},\tau_{\dio}>0$ such that
\bee
\label{eq_diophantine_definition}
\|\bn\cdot\omega\|:=\dist(\bn\cdot\omega,\Z)\ge C_{\dio}|\bn|^{-\tau_{\dio}},\quad \forall \bn\in \Z^d\setminus\{\bze\}.
\ene
Without loss of generality, we will always assume $0<\omega_1<\ldots<\omega_d<1/2$. The set of Diophantine vectors with the above property will be denoted by $\mathrm{DC}_{C_{\dio},\tau_{\dio}}$, implying that the dependence on $d$ will be clear from the context. We will only use this definition with $\tau_{\dio}>d+1$.
\subsection{Operators with convergent perturbation series}

In \cite{KPS}, it was shown that, if $f$ is $C_{\reg}$-regular on $(-1/2,1/2)$ and $\omega$ is Diophantine, then the Rayleigh--Schr\"odinger perturbation series (see below) converges for sufficiently small $\varepsilon>0$. However, one can also apply the construction from \cite{KPS} in the case where $C_{\reg}$ and $D_{\min}$ themselves depend on $\ep$, under some additional restrictions on the off-diagonal terms of $H(x)$. In \cite[Section 6]{KPS}, an example of such an operator was considered. In this section, we will describe a slightly more general class of operators for which the construction from the end of \cite[Section 6]{KPS} can be applied, virtually, without any changes. The main results of the present paper will be obtained by reducing various operators to the class described in this section.

We will consider long range quasiperiodic operators with variable hopping terms. A {\it quasiperiodic hopping matrix} is, by definition, a matrix with elements of the following form
\bee
\label{eq_quasi_1}
\Phi_{\bm\bn}(x)=\varphi_{\bm-\bn}(x+\omega\cdot(\bm+\bn)/2),\quad \bm,\bn\in \Z^d,
\ene
where $\varphi_{\bm}\colon \mathbb R\to \mathbb C$ are Lipschitz $1$-periodic functions, satisfying the self-adjointness condition:
$$
\varphi_{\bm}=\overline{\varphi_{-\bm}}.
$$
Let also
$$
\|\varphi\|_{\ep}:=\max\{\sup_{\bk} \|\varphi_{\bk}\|_{\infty},\ep\sup_{\bk} \|\varphi'_{\bk}\|_{\infty}\}.
$$
Define $\range(\Phi)$ to be the smallest number $L\ge 0$ such that $\Phi_{\bn\bm}\equiv 0$ for $|\bm-\bn|>L$. We will only consider hopping matrices of finite range. Note that \eqref{eq_quasi_1} can be reformulated as the following covariance property:
$$
\Phi_{\bm+\ba,\bn+\ba}(x)=\Phi_{\bm\bn}(x+\ba\cdot\omega),\quad \bm,\bn,\ba\in \Z^d.
$$
Fix some $R\in \mathbb N$, and suppose that $\Phi^1,\Phi^2,\ldots$ is a family of quasiperiodic hopping matrices with $\range(\Phi_k)\le k R$, defined by a family of functions $\varphi^1_{\bm}, \varphi^2_{\bm}, \ldots$. The class of operators we would like to consider will be of the following form:
\bee
\label{eq_longrange_def}
H=V+\varepsilon\Phi^1+\varepsilon^2\Phi^2+\ldots,\quad 0\le \varepsilon<1,
\ene
where
$$
(V(x)\psi)_{\bn}=f(x+\bn\cdot\omega)\psi_{\bn}.
$$
One can easily check that, assuming
$$
\|\varphi\|_{\infty}=\sup_{j}\|\varphi^j\|_{\infty}=\sup_{j,\bm}\|\varphi^j_{\bm}\|_{\infty}<+\infty,\quad 0\le\varepsilon<1,
$$
the part $\Phi=\varepsilon\Phi^1+\varepsilon^2\Phi^2+\ldots$ defines a bounded operator on $\ell^2(\Z^d)$.

The central object of \cite{KPS} is the {\it Rayleigh--Schr\"odinger perturbation series}, which is a formal series of eigenvalues and eigenvectors
\bee
\label{eq_lak}
E=E_0+\varepsilon E_1+\varepsilon^2 E_2+\ldots,
\ene
\bee
\label{eq_psik}
\psi=\psi_0+\varepsilon\psi_1+\varepsilon^2\psi_2+\ldots
\ene
where, in a small departure from the notation of \cite{KPS}, we assume
\bee
\label{eq_rs_conditions}
\lambda_0=f(x_0+\bn\cdot\omega),\quad \psi_0=e_{\bn},\quad \psi_j\perp\psi_0\,\text{ for }j\neq 0.
\ene
Under the above assumptions, we consider the eigenvalue equation
$$
H(x_0)\psi=E\psi
$$
as equality of coefficients of two power series in the variable $\ep$:
$$
(V+\varepsilon\Phi^1+\varepsilon^2\Phi^2+\ldots)(\psi_0+\varepsilon\psi_1+\varepsilon^2\psi_2+\ldots)=
(E_0+\varepsilon E_1+\varepsilon^2 E_2+\ldots)(\psi_0+\varepsilon\psi_1+\varepsilon^2\psi_2+\ldots).
$$
Assuming that $f$ is strictly monotone, the above system of equations has a unique solution satisfying \eqref{eq_rs_conditions}.

It will also be convenient to consider a graph $\Gamma(x_0)$, whose set of vertices is $\Z^d$, and there is an edge between $\bm$ and $\bn$ if $\Phi_{\bm\bn}(x_0)\neq 0$. The {\it length} of that edge is the smallest $j>0$ such that $\Phi_{\bm\bn}^j(x_0)\neq 0$.

We can now describe the class of operators with convergent perturbation series which will be used in this paper. Let $I_1,\ldots,I_k\subset (-1/2,1/2)$ be disjoint closed intervals, and $\mu_1,\ldots,\mu_k>0$. We will only consider $\mu_j\in \N$ in applications, but the argument can be easily extended for $\mu_j\in [0,+\infty)$. We will make the following assumptions.
\begin{itemize}
	\item[(conv1)] $f$ satisfies (f1) from the previous section. Additionally, $f$ is $C_{\reg}$-regular on $(0,1)\setminus (I_1\cup\ldots\cup I_k)$ with $D_{\min}\ge 1$.
	\item[(conv2)] $c_1 \varepsilon^{\mu_j}\le f'(x)\le c_2$ on $I_j$. Here, $f$ is allowed to be non-differentiable at some values of $x$, but then the inequality is required for all derivative numbers. As a consequence, $f$ is strictly monotone on $(-1/2,1/2)$.
	\item[(conv3)] Suppose that $x_0+\omega\cdot\bn\in I_j$. Any edge of $\Gamma(x_0)$ that starts at $\bn$ has length at least $\mu_j+1$.
\end{itemize}
The proof of the following result can be done along the same lines as the argument in \cite[Proof of Theorem 6.2]{KPS}.
\begin{prop}
\label{prop_convergence}
Under the above assumptions $(\mathrm{conv1})$ -- $(\mathrm{conv4})$, there exists
$$\ep_0=\ep_0(C_{\reg},\{\mu_j\}_{j=1}^k,C_{\dio},\tau_{\dio},\|\varphi\|_{\ep},c_1,c_2)>0
$$
such that, for $0<\ep<\ep_0$, the perturbation series for the operator \eqref{eq_longrange_def} is convergent, and the operator $H(x)$ satisfies Anderson localization for all $x\in \R$.
\end{prop}
\begin{rem}
Here, we use an additional refinement of \cite{KPS}, by introducing $\|\varphi\|_{\ep}$ instead of considering $\|\varphi\|_{\infty}$ and $\|\varphi'\|_{\infty}$ separately. The reason is that differentiating $\varphi$ is ``cheaper'' than differentiating the denominators, and the latter has already been accounted for in (conv2).
\end{rem}
One can also consider the perturbation series in finite volume. By ``finite volume'', we mean the case when the operator is restricted to a finite box in $\Z^d$. The results of \cite{KPS}, as well as Proposition \ref{prop_convergence}, do not extend directly into that case (in the language of \cite{KPS}, certain loops required for cancellation are forbidden). However, in all our applications, the volume will be fixed before choosing a small $\ep$, in which case one can apply the regular perturbation theory of isolated eigenvalues:
\begin{prop}
\label{prop_finite_isolated}
Let $H=V+\ep \Phi$ be a self-adjoint operator on $\ell^2([0,N]^d)$:
$$
(H\psi)_{\bn}=V_{\bn}\psi_{\bn}+\ep \sum_{\bm}\Phi_{\bn\bm}\psi_{\bm},
$$
with the diagonal part $V$ (in other words, $\Phi_{\bn\bn}\equiv 0$). Assume that all values $V_{\bk}$ are disjoint:
$$
\delta=\max_{\bm,\bk\colon\bm\neq \bk}|V_{\bm}-V_{\bk}|>0.
$$
Then the coefficients of the perturbation series \eqref{eq_lak}, \eqref{eq_psik} satisfy
$$
|E_j|\le \frac {C(N,d)^j\|\Phi\|^j}{\delta^j},\quad \|\psi_j\|\le \frac {C(N,d)^j\|\Phi\|^j}{\delta^j}.
$$
As a consequence, both series converge for small $\ep>0$.
\end{prop}
\begin{cor}
\label{cor_finite_1}
Suppose that $f$ satisfies $(f1)$, and $\{1,\omega_1,\ldots,\omega_d\}$ are rationally independent. For any $M>0$ and $C_{\reg}>0$ there exists
$$
\ep_0=\ep_0(M,C_{\reg},\omega)>0,\quad \delta=\delta(M,C_{\reg},\omega)
$$
such that the following is true. For any $x_0\in \R$, any box $B$ in $\Z^d$ with $|B|\le M$, any $\ep$ with $0<\ep<\ep_0$ and any $\bn\in B$, such that $f$ is $C_{\reg}$-regular at $x_0+\omega\cdot\bn$, the operator $H$ restricted to the box $B$ has a unique simple eigenvalue $E=E(x_0)$ with
$$
|E(x_0)-f(x_0+\omega\cdot\bn)|\le \ep \delta^{-1}\|\varphi\|_{\infty}.
$$
The corresponding eigenvector $\psi$, $\|\psi\|_{\ell^2(B)}=1$, can be chosen to have
$$
|\psi(x_0)-e_{\bn}|\le \ep \delta^{-1}\|\varphi\|_{\infty}.
$$
Moreover, both $E$ and $\psi$ are Lipschitz continuous at $x_0$ and satisfy similar bounds:
$$
|E'(x_0)-f'(x_0+\omega\cdot\bn)|\le \ep\delta^{-1}\|\varphi'\|_{\infty},\quad |\psi'(x_0)|\le \ep\delta^{-1}\|\varphi'\|_{\infty}.
$$
\end{cor}
\begin{proof}
Since the components of $\omega$ are rationally independent and $f$ is $C_{\reg}$-regular at $x_0+\omega\cdot\bn$, we have $f(x_0+\omega\cdot\bn)$ separated from the other diagonal entries by a constant that only depends on $\omega$, $C_{\reg}$, and the size of the box. Then the argument follows from Proposition \ref{prop_finite_isolated}. The differentiability follows from differentiating the perturbation series in the variable $x_0$ which, in turn, follows from (cr2).
\end{proof}
\begin{rem}
\label{rem_regular_small}
Exact same bounds on the derivatives also hold for the infinite volume eigenvectors (in that case, $\delta=\delta(C_{\dio},\tau_{\dio},d)$). As a consequence, for most practical purposes involving at most one derivative, regular eigenvalues can be considered to be small perturbations of the original diagonal entries, both in finite and infinite volumes. In a finite volume, the constants depend on the size of the box, but can be made independent of the location of the box.
\end{rem}
The following is a simple consequence of rational independence.
\begin{prop}
\label{prop_largest}
Suppose that $f$ satisfies $(f1)$. For any $M>0$ there exists $W=W(f,M,\omega)$ such that, for any box $B$ with $|B|\le M$, all eigenvalues of $H_B(x)$, except for at most one, are bounded by $W+\ep\|\varphi\|_{\infty}$ in absolute value.
\end{prop}
\begin{proof}
The only possibility for the operator to have a very large eigenvalue is for $x_0+\omega\cdot\bn$ to be close to $\Z+1/2$. However, this implies that all the remaining values of $x_0+\omega\cdot \bm$, where $\bm\neq \bn$ is from the same box, are away from $\Z+1/2$. Ultimately, the dependence on $\omega$ is only through $C_{\dio}$ and $\tau_{\dio}$.
\end{proof}
\section{Some properties of Schr\"odinger eigenvectors}

In this section, we will summarize some basic properties of the eigenvectors of discrete Schr\"odinger operators in bounded domains. Let $C\subset \Z^d$. Denote the Laplace operator on $\ell^2(C)$ (with the Dirichlet boundary conditions on $\Z^d\setminus C$) by
$$
(\Delta_C\psi)(\bn)=\sum\limits_{\bm\in C\colon |\bm-\bn|_1=1}\psi(\bm).
$$
\begin{defn}
\label{def_directly_reachable}
Let $A,B\subset \Z^d$. We will define $\bn\in A\cup B$ to be {\it directly reachable from $B$} using the following recurrent definition:
\begin{enumerate}
	\item Any point $\bn\in B$ is directly reachable in $0$ steps.
	\item Let $\bn\in A\cup B$, $\bm\in A$, $|\bm-\bn|=1$ and the set 
$$
\{\bm'\in \Z^d \colon |\bn-\bm'|\le 1,\,\bm'\neq \bm\}
$$
is contained in $A\cup B$ and only consists of directly reachable points in at most $k-1$ steps. Then we declare $\bm$ to be directly reachable in at most $k$ steps.
\end{enumerate}
\end{defn}
The above definition can be understood as follows. The eigenvalue equation 
\bee
\label{eq_eigen_equation}
\sum\limits_{\bm\colon |\bm-\bn|_1=1}\psi(\bm)=(E-V(\bn))\psi(\bn)
\ene
allows us to determine the value $\psi(\bm)$ from the values
$$
\{\psi(\bm')\colon |\bm'-\bn|\le 1,\,\bm'\neq \bm\},
$$
where $\bn\in A\cup B$ is some point with $|\bn-\bm|=1$.  A directly reachable point, therefore, is a point $\bm\in A\cup B$ such that $\psi(\bm)$ can be determined from the values of $\psi$ on $B$ using a finite sequence of such operations.

\begin{defn}
\label{def_direct_continuation}
We will say that $B\subset \Z^d$ satisfies {\it direct unique continuation property} (DUCP) for $A\subset \Z^d$ if all points of $A$ are directly reachable from $B$.	
\end{defn}


We will need the following very coarse quantitative version of the unique continuation. Better modifications to this lemma can be made in particular situations.

\begin{lem}
\label{lemma_unique}
Assume that $B\subset \Z^d$ satisfies DUCP for $A\subset \Z^d$, and any point of $A$ is directly reachable from $B$ in at most $N$ steps, and $|V_{\bn}|\le W$ for all $\bn\in A\cup B$ for some $W>0$. Assume that $H_{A\cup B}\psi=E\psi$, and $\|\psi\|_{\ell^2(A\cup B)}\ge 1$. Then there exists $\bm\in B$ such that
$$
|\psi(\bm)|\ge |A\cup B|^{-1}(2^d+|E|+W)^{-N}.
$$
\end{lem}
\begin{proof}
By contradiction: each step we can only increase the value of the eigenfunction by a factor of $2^d+|E|+W$.
\end{proof}
\begin{cor}
\label{cor_unique_epsilon}
In Lemma $\ref{lemma_unique}$, replace $\Delta$ by $\ep\Delta$. Then the conclusion can be restated as
$$
|\psi(\bm)|\ge |A\cup B|^{-1}(2^d\ep+(|E|+W))^{-N}\ep^N.
$$
\end{cor}

\begin{rem}
In Lemma $\ref{lemma_unique}$, one can generalize the Laplace operator in the following way: the weight of any edge connecting two adjacent points of $B$ can be replaced by any number between $-1$ and $1$. The reason is that, in the process of calculating the new values of $\psi_{\bm}$, $\bm\in B\setminus A$, one would not have to divide by these weights.
\end{rem}
\begin{rem}
\label{remark_example_unique}
For a subset $A\subset \Z^d$ and $k\in \mathbb N$, let
$$
\partial_k A:=\{\bn\in \Z^d\setminus A\colon \dist(\bn,A)\le k\}.
$$
One can easily check that $B:=\partial_2 A$ satisfies DUCP for any bounded $A$. Now, suppose that $A$ is a box. Then, in Corollary \ref{cor_unique_epsilon}, one can take $N=\lceil\diam(B)\rceil$. Moreover, in Corollary \ref{cor_unique_epsilon}, one can take $W$ to be the {\it second largest} element of the set
$$
\{|V_{\bn}|\colon \bn\in A\cup \partial_1 A\}.
$$
In other words, for any $\bm\in A$, one can calculate $\psi(\bm)$ from the values of $\psi$ on $B$ through a sequence of applications of \eqref{eq_eigen_equation} without using the value $V(\bn)$ of the potential at any particular point $\bn$.  Thus, the value of the potential at one point is allowed to be arbitrary large and yet excluded from consideration in the bound from \ref{cor_unique_epsilon}. 

This argument will be particularly useful in combination with Proposition \ref{prop_largest}, since it states that, while we may not be able to control the largest eigenvalue on a box of given size, we can control all remaining eigenvalues in terms of $\omega$, $f$, and the size of the box.
\end{rem}
We will also need some information about the eigenfunction decay, which follows from elementary perturbation theory.
\begin{prop}
\label{prop_clustering}
	Let $A$ and $H_A$ be as above, and let $A=A_1\sqcup\ldots \sqcup A_n$ be a disjoint union. Assume that
$$
|V_{\bm}-V_{\bn}|\ge\eta,\quad \forall \bm\in A_j,\,\,\bn\in A_k,\,\,j\neq k.
$$
Then there exists an orthonormal basis $\{\psi_{\bm}\colon \bm\in A\}$ of eigenfunctions of $H_A$ such that if $\bm\in A_j$, then
\bee
\label{eq_clusters}
|\psi_{\bm}(\bn)|\le C(A,\{A_j\})\eta^{-\dist(\bn,A_j)}.
\ene
\end{prop}
In other words, if the values of $V$ are clustered, then the corresponding eigenfunctions decay exponentially away from the clusters. The proof immediately follows from the standard perturbation theory \cite{Kato}. The statement is only meaningful for $\eta\gg 1$, otherwise one can absorb the bound into $C(A,\{A_j\})$. The constant does not depend on $V$ and on $\eta$, as long as \eqref{eq_clusters} is satisfied. In particular, it does not depend on the energy range inside a cluster; additionally, any dependence on the number of lattice points in each cluster can also be absorbed into $C$.
\begin{rem}
In Proposition \ref{prop_clustering}, one can replace the Laplacian by a weighted Laplacian as long as the weights are bounded above by $1$ in absolute value, with the same constants. One can also replace it by a long range operator if the distance function is modified accordingly.
\end{rem}
The following statement, which also follows from perturbation theory for $2\times 2$ matrices, is useful.
\begin{prop}
\label{prop_epsilon2}
Under the assumptions of Proposition $\ref{prop_clustering}$, replace $\Delta_A$ by a weighted Laplacian, with all weights bounded by $1$ in absolute value, and the weights of all edges connecting points from different clusters are bounded by $\ep$. Denote by $H_A'$ the operator $H_A$ with weights of all edges between different clusters replaced by zero. Then
$$
|\lambda_j(H_A')-\lambda_j(H_A)|\le C(A,\{A_j\}) \varepsilon^2\eta^{-1},
$$
where $\lambda_j(H)$ denotes the $j$th eigenvalue of $H$ in the increasing order.
\end{prop}
\begin{proof}
This is a standard argument from perturbation theory. For each pair of points $\bn,\bm$ from different clusters, let
$$
P=\<e_{\bn},\cdot\>e_{\bn}+\<e_{\bm},\cdot\>e_{\bm},
$$
and consider the operator 
$$
H_{\{\bm,\bn\}}:=\left.P H_A P\right|_{\ran P};
$$ 
in other words, restrict $H_A$ to $\ell^2\{e_{\bn},e_{\bm}\}$. Since $|V_{\bm}-V_{\bn}|\ge \eta$, for $\ep<\eta$ we have two eigenvectors of $H_{\{\bm,\bn\}}$ which are small perturbations of $e_{\bm},e_{\bn}$, respectively. Denote these (orthonormal) eigenvectors by $f_{\bm},f_{\bn}$, and consider the unitary matrix with columns $f_{\bm},f_{\bn}$:
$$
V=(f_{\bm},f_{\bn}).
$$
Let 
$$
V_1:=V\oplus\one_{\ell^2(A\setminus\{\bm,\bn\}},
$$ 
and
$$
H^1_A=V_1^{-1}H_A V.
$$
The above procedure is called a {\it partial $2\times 2$ diagonalisation}. One can check that, after such transformation, the edge between $\bm$ and $\bn$ of length $\ep$ would disappear, and all other entries of the operator are changed at most by $O(\ep^2\eta^{-1})$ (it is possible that the edge between $\bm$ and $\bn$ does not completely disappear, but its new length will be at most $O(\ep^2\eta^{-1})$).

Apply partial $2\times 2$ diagonalisation to each edge that connects a pair of vertices from different clusters. By doing this, one will reduce all weights of those edges to at most $O(\ep^2\eta^{-1})$. Afterwards, and additional perturbation of the whole operator of size $O(\ep^2\eta^{-1})$ would completely eliminate those edges. Finally, one can restore the entries inside of each cluster to their original form $H_A'$ by another perturbation of the order $O(\ep^2\eta^{-1})$. Note that the constants do depend on the size of the domain and geometry of clusters; however, they can be chosen uniformly in the magnitude of the potential (in fact, large values of the potential only make the situation easier), as long as the clusters are $\eta$-separated.
\end{proof}
We will also need the following property related to the differentiation of eigenvalues with respect to a parameter. The following is often referred to as Hellmann--Feynman variational formula, see \cite[Remark II.2.2]{Kato}
\begin{prop}
\label{prop_hell}
Let $t\mapsto A(t)$ be a Lipschitz continuous family of self-adjoint matrices. Let $\lambda(t)$ be a branch of an isolated eigenvalue of $A(t)$, and let $\psi(t)$ be a corresponding $\ell^2$-normalized eigenvector. Then
\bee
\label{eq_hell}
\lambda'(t)=\frac{d}{dt}\<A(t)\psi(t),\psi(t)\>=\<A'(t)\psi(t),\psi(t)\>.
\ene
\end{prop}
\begin{proof}
The remaining term in the right hand side of \eqref{eq_hell} is $\lambda(t)\frac{d}{dt}\langle \psi(t),\psi(t)\rangle $, which is identically zero due to normalization of $\psi(t)$.
\end{proof}
In particular, if $A(t)=A_0+f(t)\<e_j,\cdot\>e_j$, then
$$
\lambda'(t)=f'(t)|\langle \psi(t),\be_j\rangle|^2.
$$

\begin{prop}
\label{prop_epsilon_loss}	
Under the assumptions of Proposition $\ref{prop_clustering}$, assume 
$$
H(t)=H_0+f(t)\<e_{\bk},\cdot\>e_{\bk}.
$$
Suppose that $E=E(t)$ is an isolated eigenvalue of $H(t)$ that is separated from the rest of the spectrum by $\delta>0$. Assume that $\bm\in A_r$, $\bn\in A_{s}$. Denote by $P(t)$ the spectral projection onto the eigenspace associated with the eigenvalue $E(t)$. Then
$$
\left|\frac{dP_{\bm\bn}(t)}{dt}\right|\le C(A,\{A_j\}) |f'(t)|\delta^{-1}\eta^{-\dist(A_r,\{\bk\})-\dist(\{\bk\},A_{s})}.
$$
\end{prop}
\begin{proof}
The spectral theorem implies
\bee
\label{eq_resolvent1}
|\<(H(t)-E)^{-1}e_{\bm},e_{\bk}\>|\le C \delta^{-1}\eta^{-\dist(A_r
,\{\bk\})}.
\ene
We also have
$$
\frac{d}{dt}P(t)_{\bm\bn}=\left(\frac{d}{dt}\frac{1}{2\pi i}\oint_{|\lambda-E|=\delta/2}(H_A-\lambda)^{-1}d\lambda\right)_{\bm\bn},
$$
from which one can get the inequality by differentiating.
\end{proof}
\begin{rem}
Our use of discrete unique continuation is relatively elementary. We refer the reader to \cite{BK,DS} for more advanced applications of unique continuation for discrete Schr\"odinger operators.
\end{rem}
\section{Matrix functions with Maryland-type diagonal entries}
Let $f\colon (-1/2,1/2)\to \R$ satisfy (f1) from Section 2. We will say that $f$ is {\it locally Lipschitz monotone at }$x_0\in (-1/2,1/2)$ if 
$$
f'(x)\ge 1,\quad \forall x\colon |f(x)-f(x_0)|\le 1.
$$
This condition is a weaker version of $C_{\reg}$-regularity at $x_0$ with $D_{\min}\ge 1$ and will be sufficient for the construction in the current section. Let also
\bee
\label{eq_homotopy}
f_t(x)=f(x)+t\{x-1/2\},
\ene
where $\{\cdot\}$ denotes the fractional part. Assume that $f,f_t$ are extended into $\R$ by $1$-periodicity. It is easy to see that $f_1$ is locally Lipschitz monotone at all $x_0\in (-1/2,1/2)$. Moreover, if $f_t$ is locally Lipschitz monotone at $x_0$, then all $f_{t'}$, $t'\ge t$, are also locally Lipschitz monotone at $x_0$.

In the following theorem, we consider $N\times N$ matrix valued functions whose diagonal entries are locally Lipschitz monotone. Let $\{e_1,\ldots,e_{N}\}$ denote the standard basis in $\C^N$.
\begin{thm}
\label{th_diagonalisation1}
Let 
$$
d_1,\ldots,d_N\in (-1/2,1/2),\quad \min_{i\neq j}|d_i-d_j|=\delta>0.
$$
Let also
$$
A_t(x)=\diag\{f_t(x-d_1),\ldots,f_t(x-d_n)\}+\ep \Phi(x),\quad A(x)=A_0(x),
$$
where $\Phi(\cdot)$ is a continuous real symmetric $1$-periodic $N\times N$ matrix-valued function on $\R$. Assume that all eigenvalues of $A_t$ are simple for all $t\in [0,1]$ and all $x\in \R$. Then there exist $c_j=c_j(N)>0$, $j=1,2$ such that, for
$$
0\le \ep\|\Phi\|_{\infty}<c_1(N)\delta,
$$
we have the following:
\begin{enumerate}
	\item Suppose that $f$ is locally Lipschitz monotone at $x_0-d_j$. Then $A(x)$ has a unique eigenvalue $E_j(x)$ and a unique $\ell^2$-normalized real eigenvector $\psi_j(x)$ satisfying:
\bee
\label{eq_middle_perturbative}
|E_j(x)-f(x-d_j)|\le c_2(N) \ep \|\Phi\|_{\infty}\delta^{-1},\quad |\psi_j(x)-e_j|\le c_2(N) \ep \|\Phi\|_{\infty}\delta^{-1}.
\ene
\item There exists a continuous $1$-periodic real orthogonal matrix function $U(\cdot)\colon \R\to O(N,\R)$ such that $U^{-1}(x)A(x)U(x)$ is a diagonal matrix and $j$th column of $U(x)$ coincides with $\psi_j(x)$ from $(1)$ for all $x$ where $(1)$ is applicable.
\end{enumerate}
\end{thm}
\begin{rem}
One needs to specify the exact meaning of ``simple'' in the case when one of the matrix entries becomes infinite. It is easy to see that, as $x$ approaches $d_j$, exactly one eigenvalue of $A_t(x)$ approaches $\pm \infty$, and the remaining eigenvalues remain bounded. In fact, they approach the eigenvalues of $A_t$ with $j$th row and column removed, and thus can be extended continuously through $d_j$. Thus, we can require for all those finite eigenvalues to be simple. Equivalently, we can require the distances between eigenvalues of $A_t(x)$ to be bounded by some constant $\varkappa>0$ for all $t\in [0,1]$ and all $x\in \R\setminus(\{d_1,\ldots,d_N\}+\Z)$.
\end{rem}
\begin{proof}
Part (1) follows from the standard perturbation theory of isolated eigenvalues: see also Proposition \ref{prop_finite_isolated}. Note that, once we fix a real $\ell^2$-normalized branch, it must be close to $e_j$ or $-e_j$, and requiring it to be close to $e_j$ completely determines this branch.

Since $f_1(\cdot)$ is locally Lipschitz monotone on $(-1/2,1/2)$, there exist (unique real normalized) continuous branches of eigenvectors:
$$
\psi_j(\cdot)\colon (-1/2+d_j,1/2+d_j)\to \R^d,\quad \|\psi_j(x)-e_j\|\le c_3(N)\ep \|\Phi\|_{\infty}\delta^{-1}.
$$
Since the spectra of $A_t$ are simple, each $\psi_j$ has a unique real $\ell^2$-normalized extension in the parameter $t$ into $[0,1]$. Denote this extension by $\psi_j(x,t)$, and the corresponding eigenvalue branch by $E_j(x,t)$. For any $x_0$ such that $f$ is locally Lipschitz monotone at $x_0-d_j$, any real eigenvector of $A_t(x_0)$ with eigenvalue $E_j(x_0,t)$ must be close to $e_j$ or $-e_j$, uniformly in the parameter $t$. Since it is close to $e_j$ for $t=1$, the same holds for all $t$.

A similar argument can be applied to $x_0$ near the endpoints. Since one (and only one) eigenvalue of $A_t(x)$ approaches infinity, the corresponding branch of eigenvector must approach $e_j$ or $-e_j$ as $x\to \pm 1/2+d_j$, uniformly in $t$. Since for $t=1$ it approaches $+e_j$ from both sides, the same holds for all $t$. In other words, it cannot suddenly change from a vector close to $e_j$ to a vector close to $-e_j$, due to continuity in $t$ and the fact that it has to be close to $e_j$ or $-e_j$ for any particular $t$.
\end{proof}
\begin{rem}
\label{remark_t_infinity}
Under the assumptions of Theorem \ref{th_diagonalisation1}, one can take $t\to +\infty$ in \eqref{eq_homotopy}. In this case, $\psi_j(x,t)$ will be well defined for all $t\ge 0$ and uniformly approach $e_j$ as $t\to\infty$. As a consequence, if the family $\{A_t(x)\}$ satisfies the assumptions of Theorem \ref{th_diagonalisation1}, there is a unique family $\{U_t(x)\}$ of real orthogonal matrices such that $U_t(x)$ diagonalises $A_t(x)$ and $U_t(x)\to I$ as $t\to +\infty$. In this notation, we will have $U(x)=U_0(x)$. One can use this as an alternative definition of $U(x)$.
\end{rem}
\begin{prop}
\label{prop_jacobi_separation} Fix $L>0$, $a\in \R$. Let $J$ be an $N\times N$ matrix satisfying the following properties: $J_{ij}=0$ for $|i-j|>1$, and $J_{ij}=J_{ji}\in \R$. Assume that $J_{i,i+1}\ge 1$ for all applicable values of $i$, and that $J_{ii}\in [a,a+L]$ for all $i$. Then the spectrum of $J$ is simple. Moreover, if $E_i\neq E_j$ are two eigenvalues of $J$, then $|E_i-E_j|\ge c(N,L)>0$.
\end{prop}
\begin{proof}
The fact that the spectrum is simple is standard and follows from the fact that the any solution to the eigenvalue equation is uniquely determined from its value at one point. The rest follows from compactness.
\end{proof}
\begin{rem}
By rescaling, one can obtain a version of the proposition with 	$|E_i-E_j|\ge c(N,\ep L)\varepsilon$ provided that $|J_{i,i+1}|>\varepsilon>0$.
\end{rem}
\begin{rem}
\label{rem_half_bounded}
In Proposition \ref{prop_jacobi_separation}, one can relax the requirement $J_{ii}\in [a,a+L]$ for the endpoint values $J_{11}$ and $J_{NN}$: one can remove this restriction for one of them, or for both if one of them goes to $+\infty$, and the other to $-\infty$.
\end{rem}
\begin{rem}
\label{rem_quasi_1d}
In higher dimensions, Proposition \ref{prop_jacobi_separation} will be often used in combination with Proposition \ref{prop_epsilon2}: if a cluster from the latter proposition has a linear shape (or is somewhat one-dimensional), then Proposition \ref{prop_jacobi_separation} guarantees that its eigenvalues will be $c\ep$-separated (after restricting the operator to the cluster). The presence of coupling terms with  other clusters will shift the eigenvalues by at most $O(\ep^2)$, and will thus preserve the separation condition.
\end{rem}
\section{The moving block construction: general scheme}
As a reminder, the standard basis of $\R^d$ (and $\Z^d$) is denoted by $\{\be_1,\be_2,\ldots,\be_d\}$. The translation operator on functions is defined by
$$
T^{\bn}\colon \ell^2(\Z^d)\to \ell^2(\Z^d),\quad (T^{\bn}\psi)(\bm)=\psi(\bm+\bn).
$$
Since we will use the following specific translations more often, we also denote
$$
T:=T^{\be_1},\,\,T^{*}:=T^{-1}=T^{-\be_1}.
$$

Consider the operator \eqref{eq_h_def} on $\ell^2(\Z^d)$ and assume that $f$ satisfies (f1). Recall that, without loss of generality, we assumed
$$
0<\omega_1<\omega_2<\ldots<\omega_d<1/2.
$$
Fix some $x_0\in \R$ and $C_{\reg}>0$. The construction below will depend on the choice of $x_0$; however, any $x_0$ would be equally applicable. A lattice point $\bn\in \Z^d$ will be called {\it regular for $H(x)$} if $f$ is $C_{\reg}$-regular at $x+\omega\cdot\bn$. As discussed in Corollary \ref{cor_finite_1} and Remark \ref{rem_regular_small}, for sufficiently small $\ep>0$, $H(x)$ has an eigenvalue close to $f(x+\omega\cdot\bn)$ and the corresponding eigenvector close to $e_{\bn}$ at any regular point $\bn$. The same also holds for a restriction of $H$ onto any box containing $\bn$. The smallness of $\ep$ depends on the size of the box. Additionally, ``close'' in this case also applies to derivatives in $x$. A lattice point which is not regular for $H(x)$ will be, naturally, called {\it singular for $H(x)$}. Let
$$
S_{\sing}(x):=\{\bn\in \Z^d\colon\bn\text{ is singular for }H(x)\}.
$$
For the purposes of the construction, the point with infinite potential will be considered regular. Since the components of $\omega$ are rationally independent, there is at most one such point. We will always assume that singular points are contained in a finite range of energies:
\begin{itemize}
	\item[(gen0)] There is $E_{\reg}>0$ such that $f$ is $C_{\reg}$-regular at any $x$ with $|f(x)|\ge E_{\reg}$.
\end{itemize}
Let $S\subset \Z^d$ be a finite subset possibly containing some singular points. Let also $R\supset S$ be another finite subset, with the following properties:
\begin{itemize}
	\item[(gen1)] $S\cup(S+\be_{1})\subset R\cap (R+\be_1)$.
	\item[(gen2)] For every $x\in [x_0,x_0+\omega_1]$, all points of $(R\cup (R+\be_1))\setminus(S\cup (S+\be_{1}))$ are regular for $H(x)$.
\end{itemize}
Let also 
\bee
\label{eq_R_def}
R_-:=R\setminus(R+\be_1),\,\, R_+:=(R+\be_1)\setminus R,\,\, R_0:=R\cap(R+\be_1),\,\, R':=R\cup(R+\be_1).
\ene 
We can represent $R'$ as a disjoint union $R'=R_-\sqcup R_+\sqcup R_0$. Denote by $H'_{R'}(x)$ the following family of operators on $\ell^2(R')$:
$$
H'_{R'}(x_0):=V_{R_-}(x_0)\oplus H_{R_0\cup R_+}(x_0);
$$
$$
H'_{R'}(x_0+\omega_1):=H_{R_-\cup R_0}(x_0+\omega_1)\oplus V_{R_+}(x_0+\omega_1),
$$
where $V_{R_{\pm}}(x)$ denotes just the potential part of $H_{R_{\pm}}(x)$, without the Laplacian, restricted onto the domains $R_{\pm}$. For the values $x\in [x_0,x_0+\omega_1]$, we will linearly interpolate:
$$
H'_{R'}(x):=(1-s)\left(V_{R_-}(x_0)\oplus H_{R_0\cup R_+}(x_0)\right)+s\left(H_{R_-\cup R_0}(x_0+\omega_1)\oplus V_{R_+}(x_0+\omega_1)\right),
$$
$$
x=x_0+s\omega_1 \quad \text{for}\,\,s\in [0,1].
$$
Assume the following condition:
\begin{itemize}
	\item[(gen3)] Matrix families $H'_{R'}(x)$, $H_{R_0\cup R_+}(x)$, $H_{R_-\cup R_0}(x)$	satisfy the assumptions of Theorem \ref{th_diagonalisation1}, for all $x\in [x_0,x_0+\omega_1]$.
\end{itemize}
Technically, in order to apply Theorem \ref{th_diagonalisation1}, $H'_{R'}(x)$ has to be defined (possibly with infinite entries) and be periodic for all $x\in \R$. The exact form of the extension will not matter. For the sake of clarity, assume that, for $x\in [x_0+\omega_1,x_0+1]$ the  off-diagonal parts of $H'_{R'}(x)$ are linearly interpolated between those of $H'_{R'}(x_0+\omega_1)$ and $H'_{R'}(x_0)$, so that $H'_{R'}(x_0+1)=H'_{R'}(x_0)$, and then extended $1$-periodically to $\R$. The later construction (see the definition of the operator $U_0^{\mb}(x)$ below) we will only use $H'_{R'}(x)$ for $[x_0,x_0+\omega_1]$.

From the definitions, we have
\bee
\label{eq_translation_relation1}
H_{R_0\cup R_+}(x_0)=T^*(H_{R_-\cup R_0}(x_0+\omega_1))T,
\ene
where $T\colon \ell^2(R_0\cup R_+)\to\ell^2(R_-\cup R_0)$ is the restriction of the translation map (which is a bijection between $R_0\cup R_+$ and $R_-\cup R_0$). Denote by $U_{R'}(x)$, $x\in [x_0,x_0+\omega_1]$, the result of applying Theorem \ref{th_diagonalisation1} to the family $H'_{R'}(x)$. For $x=x_0$ and $x=x_0+\omega$, $U_{R'}(x)$ splits into a direct sum:
$$
U_{R'}(x_0)=\one_{R_-}\oplus U_{R_0\cup R_+}(x_0);
$$
$$
U_{R'}(x_0+\omega_1)=U_{R_-\cup R_0}(x_0+\omega_1)\oplus \one_{R_+}.
$$
The matrices $U_{R_0\cup R_+}(x_0)$ and $U_{R_-\cup R_0}(x_0+\omega_1)$ have real entries and, due to \eqref{eq_translation_relation1}, diagonalise the same operator, up to the identification of $R_-\cup R_0$ with $R_0\cup R_+$. Therefore, they consist of the same eigenvectors, up to the ordering and the choice of the signs. On the other hand, (gen3) guarantees that both approach the identity as $t\to +\infty$, see Remark \ref{remark_t_infinity}. Therefore, they must be related in the same way as the operators \eqref{eq_translation_relation1} they diagonalise:
\bee
\label{eq_translation_relation2}
U_{R_0\cup R_+}(x_0)=T^*U_{R_-\cup R_0}(x_0+\omega_1)T.
\ene
Extend $U_{R'}(x)$ into $\ell^2(\Z^d\setminus R')$ by identity:
\bee
\label{eq_support}
U_0^{\mathrm{mb}}(x):=U_{R'}(x)\oplus\one_{\Z^d\setminus R'}.
\ene
Then we have
\bee
\label{eq_covariance1}
T^*U_{0}^{\mb}(x_0+\omega_1)T=U_{0}^{\mb}(x_0),
\ene
We start from $U_{0}^{\mb}(x)$ defined for $x\in [x_0,x_0+\omega_1]$. One can easily check that there is a unique extension of this family into $x\in \R$ satisfying 
\bee
\label{eq_covariance2}
U_{0}^{\mb}(x+\omega_1)=T^{-1} U_{0}^{\mb}(x) T,\quad \forall x\in \R.
\ene
This above construction will be called {\it the diagonalisation of a moving block}, which is reflected in the superscript $\mb$. Note that the singular set also has the following covariance property:
\bee
\label{eq_covariance3}
S_{\sing}(x_0)=S_{\sing}(x_0-\omega\cdot\bn)+\bn.
\ene
As a consequence, if the above-mentioned set $S$ from the moving block construction contains singular points for $H(x_0)$, they will naturally correspond to singular points in $S-\be_1$ for $H(x_0+\omega_1)$. In the above construction, as $x$ increases, the set $S$ `moves' along $\Z^d$ with speed $1/\omega_1$ and traces a part of $S_{\sing}$, thus justifying the name of the construction.

It will also be convenient to use the following language: a normal operator on $\ell^2(\Z^d)$ is {\it supported} on $R'$ if $\ell^2(R')$ is its invariant subspace and the operator acts as an identity on $\ell^2(\Z^d\setminus R')$. The above construction implies that $U_{0}^{\mb}(x)$ is supported on $R'$ for $x\in [x_0,x_0+\omega_1]$ and is supported on $R'-m\be_1$ for $x\in [x_0+m\omega_1,x_0+(m+1)\omega_1]$, where $m\in \Z$.

The operator family $U_0^{\mb}(x)$ takes care only of a part of $S_{\sing}(x)$ that corresponds to a finite block $S$ that moves as $x$ moves around $\R$. We now need to cover the whole singular set by such blocks, and multiply the corresponding operators. At the same time, the family $U_{0}^{\mb}(x)$ is not a quasiperiodic operator family: its entries are not even periodic in $x$. We will need to modify this family in order for it to be quasiperiodic. First, note that the set $S_{\sing}(x)$ is not only translationally covariant \eqref{eq_covariance3}, but also translationally invariant:
$$
S_{\sing}(x+1)=S_{\sing}(x).
$$
Therefore, with each moving block, we actually have a ``train'' of blocks, with a translation by $1$. Thus, it is natural to consider
\bee
\label{eq_forced_periodic}
U^{\mb}_1(x):=\prod_{m\in \Z}U_{0}^{\mb}(x+m).
\ene
We will later require (see (gen4) below) that the supports of $U_{0}^{\mb}(x+m)$ are all disjoint. As a consequence, the above product will be well defined.

Unlike $U_{0}^{\mb}(x)$, the family $U_{1}^{\mb}(x)$ is actually a periodic function in $x$, and is quasiperiodic with respect to translations by $\omega_1$:
\bee
U_{1}^{\mb}(x)=U_1^{\mb}(x+1),\quad U^{\mb}_{1}(x+\omega_1)=T^{-1} U^{\mb}_{1}(x) T,\quad \forall x\in \R.
\ene
We now enforce the full covariance (quasiperiodicity) condition. Let $\omega'=(\omega_2,\ldots,\omega_d)$. For any $\bn'\in\Z^{d-1}\setminus\{\bze\}$, the family $x\mapsto T^{(0,\bn')}V_1(x+\bn'\cdot\omega')T^{-(0,\bn')}$ describes a different moving block (in fact, a different train of moving blocks described above). Let
\bee
\label{eq_covariant_block}
U^{\mb}_2(x):=\prod_{\bn'\in\Z^{d-1}} T^{(0,\bn')}U_1(x+\bn'\cdot\omega')T^{-(0,\bn')}.
\ene
Each factor in \eqref{eq_covariant_block} contains countably many copies of $U^{\mb}_0(x)$. If well defined, $U^{\mb}_2(x)$ is a fully quasiperiodic operator family:
\bee
\label{eq_u_quasiperiodic}
U^{\mb}_2(x)=U^{\mb}_2(x+1),\quad U^{\mb}_{2}(x+\omega\cdot\bn)=T^{-\bn} U^{\mb}_{2}(x) T^{\bn},\quad \forall x\in \R.
\ene
If order for $U_1^{\mb}$ and $U^{\mb}_2$ to be well defined, we will require the last condition:
\begin{itemize}
	\item[(gen4)]The supports of each copy of $U_0^{\mb}$ inside of the operator \eqref{eq_covariant_block} and \eqref{eq_forced_periodic} are disjoint.
\end{itemize}
The above construction will be referred to as a {\it covariant family of moving blocks}. The set $S$ will be called the {\it base block}. 
\begin{rem}
\label{remark_support_location}
Once we fix $S$ and $R'$ for $x\in [x_0,x_0+\omega_1]$, the position of other copies of $S$ and $R'$, appearing in different factors of \eqref{eq_covariant_block} and \eqref{eq_forced_periodic}, is completely determined by $S$, $R$, and $\omega$.
\end{rem}

One may need more than one such family to take care of all $S_{\sing}(x)$, in which case we will require all involved operators $U_2^{\mb}$ to commute with each other (in other words, all sets $R'$ under consideration will not overlap with each other). For simplicity, we will now consider the case with only one covariant moving block. Let
$$
H_2(x):=U^{\mb}_2(x)^{-1}H(x)U^{\mb}_2(x).
$$
We will now outline the general strategy of the proof of localization. Our goal is to show that, under some additional conditions on the sets $R_0$ (which will impose some separation conditions on $S$), the operator family $H_2(x)$ will satisfy the assumptions of Proposition \ref{prop_convergence}. For the convenience of referencing, we will denote the steps of the strategy by (strat1)--(strat5).
\begin{itemize}
	\item[(strat1)] Quasiperiodicity of $U_2^{\mb}(x)$ \eqref{eq_u_quasiperiodic} implies that $H_2(x)$ is also a quasiperiodic operator family, with diagonal entries defined by a new function $f_2$:
	$$
	(H_2(x))_{\bn\bn}=f_2(x+\omega\cdot\bn).
	$$
	For sufficiently small $\ep$ (depending on the size of $R'$), the function $f_2$ will be, say, $C_{\reg}/2$-regular at any $x$ where the original $f$ was $C_{\reg}$-regular, see Corollary \ref{cor_finite_1} and Remark \ref{rem_regular_small}.
	\item[(strat2)] As a consequence, our attention should be towards the diagonal entries of $H_2$ corresponding to singular points for $H(x)$. Assume that $\dist((S\cup (S+\be_1)),\Z^d\setminus R_0)\ge 1$, so that we avoid the boundary effects. Then, for any $\bm\in (S\cup (S+\be_1))$, $f_2(x+\omega\cdot\bm)$ is the exact eigenvalue of $H'_{R'}(x)$ that corresponds to the lattice point $\bm$. Here, we use the relation between eigenvalues of $H'_{R'}(x)$ and lattice points of $R'$ provided in Theorem \ref{th_diagonalisation1}.
	
	\item[(strat3)] Suppose that, for a single block, we have a subset $B\subset R_0\setminus(S\cup (S+\be_1))$ that satisfies DUCP for $S\cup (S+\be_1)$. For example, if $R_0$ is a box whose boundary is not too close to $S\cup (S+\be_1)$, then we can surround $S\cup (S+\be_1)$ by a layer of thickness two inside $R_0$, see Remark \ref{remark_example_unique}. For any eigenvalue of $H'_{R'}(x)$, corresponding to a lattice point $\bm \in S\cup (S+\be_1)$, apply Lemma \ref{lemma_unique} to the corresponding eigenfunction, and conclude that it must be non-trivially supported on $B$. Since every point of $B$ is regular, the corresponding diagonal entry of the operator is Lipschitz monotone. Proposition \ref{prop_hell} guarantees that there is a non-trivial positive contribution to the derivative of each singular eigenvalue from some diagonal entries of $B$. The contribution from the remaining diagonal entries of $H'_{R'}(x)$ to the derivative is also non-negative.
	\item[(strat4)] The previous part almost guarantees (conv2), with one caveat: the interpolation process in the definition of $H'_{R'}(x)$ creates some $x$-dependent off-diagonal entries, which also contribute to the derivatives of singular eigenvalues not necessarily in a positive way. However, these entries are outside of $R_0$. If we require $S\cup (S+\be_1)$ to be away from $\Z^d\setminus R_0$, we can make this contribution arbitrarily small, independent of the lower bound on the contribution from $B$ described above. This will complete the verification of (conv2). The particular selection of the set $B$ does not affect the definition of $U_2^{\mb}$, it only affects our choice of what to use as a lower bound on the derivatives. Thus, it can be chosen in an $x$-dependent way, as long as the ultimate lower bounds are uniform in $x$.
	\item[(strat5)] In order to verify (conv3), note that, after diagonalisation, singular points will only be coupled with points outside of $R_0$. The strength of this coupling (i.e, in the language of Section 2, the length of the corresponding edges) will be in terms of the decay of the corresponding eigenfunctions of $H'_{R'}$ which, in turn, can be obtained using Proposition \ref{prop_clustering} and Proposition \ref{prop_epsilon_loss} (the latter is used to control the derivatives, as explained in (conv4)). Ultimately, one can make it bounded by an arbitrarily large power of $\ep$, by requiring $\dist((S\cup (S+\be_1)),\Z^d\setminus R_0)$ to be large enough.
\end{itemize}

We will summarize the above construction in the following theorem.

\begin{thm}
\label{th_main_abstract}
Fix $d$, $C_{\reg}$, $C_{\dio}$, $\tau_{\dio}$, $D_{\min}\ge 1$, $c_{\mathrm{sep}}>0$. Let also $\omega\in \mathrm{DC}_{C_{\dio},\tau_{\dio}}$ satisfy 
$$
0<\omega_1<\omega_2<\ldots<\omega_d<1/2.
$$
Let $S\subset \Z^d$ be a finite subset. There exists $r=r(S)\in \N$ and
$$
\ep_0=\ep_0(d,C_{\reg},C_{\dio},\tau_{\dio},c_{\mathrm{sep}},S,E_{\reg},f)>0
$$
such that, for $0<\ep<\ep_0$ and all $x\in \R$, the operator $H(x)$ satisfies Anderson localization and is unitarily equivalent to an operator with convergent perturbation series, if the following are satisfied for some $x_0\in \R$:
\begin{enumerate}
	\item $f$ satisfies $(f1)$ and $\mathrm{(gen0)}$.
	\item There exists a box $R\supset S$ such that the family of covariant moving blocks, constructed above, satisfies $\mathrm{(gen1)}$--$\mathrm{(gen4)}$.
	\item Define $R_0$, $R'$ as in $\eqref{eq_R_def}$. We will require $\dist(S\cup (S+\be_1),\Z^d\setminus R_{0})\ge r$.
	\item $S_{\sing}(x)$ is contained in the support of $U_2^{\mb}(x)$, constructed above in $\eqref{eq_covariant_block}$. Due to $\mathrm{(gen2)}$, this means that $S_{\sing}(x)$ is contained in the union of the copies of the sets $S\cup(S+\be_1)$ that appear in the definition of $U_2^{\mb}(x)$.
	\item The eigenvalues of $H_{S\cup (S+\be_1)}(x)$ are uniformly $c_{\mathrm{sep}}\ep$-separated for all $x\in [x_0,x_0+\omega_1]$. Moreover, this property holds with $f$ replaced by $f_t$ from \eqref{eq_homotopy}, uniformly in $t\in [0,1]$.
\end{enumerate}
\end{thm}
\begin{rem}
\label{rem_f_dependence} There is additional dependence of $\ep_0$ on $f$, besides the dependence through $C_{\reg}$ and $E_{\reg}$. This dependence can be made explicit; we will specify it here in order to avoid overloading the statement of the main result. Consider the smallest box $B$ such that $\dist(S\cup (S+\be_1),\Z^d\setminus B)\ge 4$. Then, $\ep_0$ depends on $f$ through the quantity $W(f,M,\omega)$ from Proposition \ref{prop_largest}. Note that the dependence on $\omega$ can be reduced to dependence on $C_{\dio},\tau_{\dio}$, and the dependence on $M$ is accounted for in the dependence on $S$. See also Remark \ref{remark_example_unique}.
\end{rem}
\begin{proof}
We will show that the procedure described in the beginning of the section results in an operator unitarily equivalent to $H(x)$ and satisfying the assumptions of Proposition \ref{prop_convergence}. Fix $r\ge 4$, and take some box $R$ such that $\dist(S\cup (S+\be_1),\Z^d\setminus R')\ge r$, where $R'$ is defined in \eqref{eq_R_def}. We will specify additional conditions on $r$ later. Fix some $x_0$ and construct the covariant family of moving blocks, as described earlier in this section. For $x\in [x_0,x_0+\omega_1]$, let $B(x)$ be the smallest box in $\Z^d$ such that
\begin{enumerate}
	\item $B(x)\supset (S\cup (S+\be_1))$, and $\dist((S\cup (S+\be_1)),\Z^d\setminus B(x))\ge 1$.
	\item For any box $R\supset S$ with $\dist(\Z^d\setminus R,S\cup (S+\be_1))\ge 4$, the set 
$$
\{\bn\in B(x)\colon \dist(\bn,\Z^d\setminus B(x))\le 1\}
$$ 
does not contain the point of $R'$ with largest absolute value of the potential. 
\end{enumerate}
The choice of $B(x)$ can be made independently of $R$. Indeed, we can always choose two non-overlapping layers of thickness two around $S$, and then pick one of them with the smaller maximal value of the potential. For example, let $B_1(x)$ be the smallest box such that $\dist((S\cup (S+\be_1)),\Z^d\setminus B_1(x))\ge 2$, and $B_2(x)$ be the smallest box such that $\dist(B_1(x),\Z^d\setminus B_2(x))\ge 2$. Then, for any $x$, we can either pick $B=B_1(x)$ or $B=B_2(x)$.

Construct the operator
$$
H_2(x)=U^{\mb}_2(x)^{-1}H(x) U^{\mb}_2(x)
$$
as earlier in the section, using $S$ and $R'$. The diagonal entries of $H_2(x)$ correspond to the lattice points and are ordered in the same way as the values of the potential. The entries at distance at least 2 from $\Z^d\setminus R'$ are exact eigenvalues of the block $H_{R'}(x)$ (for the respective lattice points). The corresponding eigenvectors have cluster structure, in the sense of Proposition \ref{prop_clustering}:  one cluster would be $(S\cup (S+\be_1))\cap S_{\sing}(x)$, and each of the remaining lattice points of $R'$ corresponds to a separate cluster, with separation bounded below by a quantity that, ultimately, depends only on $r$ but not on $\varepsilon$. Proposition \ref{prop_clustering} implies that each regular eigenvector $\psi$ decays as
$$
\psi(\bm)\le C \varepsilon^{\dist(\bm,\bn)},
$$
where $\bn$ is its support point, and each singular eigenvector decays as
$$
\psi(\bm)\le C \varepsilon^{\dist(\bm,S\cup(S+\be_1))}.
$$
Here, $C$ depends on $r$, but not on $\ep$. Moreover, the derivatives of the eigenvectors in $x$ satisfy similar decay with a loss of one $\varepsilon$, see Proposition \ref{prop_epsilon_loss}.

As a consequence, each singular diagonal entry of $H_2(x)$ located on $S\cup (S+\be_1)$ has only outgoing edges of length bounded by $C\varepsilon^{r-2}$, with a derivative bound $C\ep^{r-3}$. It would be sufficient to show that the assumption (conv2) is satisfied with $\mu$ independent of $r$. Note that we already have $\mu=1$ for each regular eigenvalue.

In order to establish the above, apply Lemma \ref{lemma_unique} with
$$
B=\{\bn\in B(x)\colon \dist(\bn,\Z^d\setminus B(x)\le 2\}, \quad A=S\cup (S+\be_1).
$$ 
Proposition \ref{prop_clustering} implies that any $\ell^2(R')$-normalized eigenfunction $\psi$ associated to a singular lattice point, satisfies 
$$
\|\psi\|_{\ell^2((S\cup (S+\be_1))\cap S_{\sing}(x)}\ge 1-O(\ep)\ge 1/2
$$ 
for sufficiently small $\ep$. Lemma \ref{lemma_unique}, in view of Remark \ref{remark_example_unique} (see also Remark \ref{rem_f_dependence}), implies that $|\psi(\mathbf{b})|\ge c_1\ep^{c_2}$ for some $\mathbf{b}\in B$, where $c_2$ only depends on $S$ and $c_1$ depends on $E_{\reg}$. Since $\mathbf{b}$ is a regular point, Proposition \ref{prop_hell} implies that the derivative of the singular eigenvalue, corresponding to $\psi$, has a contribution bounded below by $C_1 \ep^{2c}$, where $C_1>0$ can be chosen to only be dependent on $S$, $\omega$, and $E_{\reg}$. The contributions from all other diagonal entries of $H'(x)$ are also non-negative. Following (strat4), we also need to account for the terms associated to the interpolated diagonal entries. These entries are located somewhere in $R_+\cup R_-$, which are outside of $R_0$ and therefore away from $S$. By combining Propositions \ref{prop_clustering} and \ref{prop_hell}, we can see that their contribution to the derivatives of singular eigenvalues is bounded above by a power of $\ep$ that grows with $r$. Thus, as long as, say, $r>3c+2$, the assumptions are satisfied.
\end{proof}
One can naturally extend the result of Theorem \ref{th_main_abstract} to the case where one requires more than one type of blocks to cover $S_{\sing}$, assuming that they do not ``interact'' with each other.
\begin{cor}
Let 
$$
\mathcal S=\{S_1,\ldots,S_k\}
$$
be a collection of subsets of $\Z^d$. There exists $r=r(\mathcal S)\in \N$ such that the following is true. Suppose that the assumptions of Theorem $\ref{th_main_abstract}$ are satisfied for each $S_j$, with a modification of $(4)$ that $S_{\sing}$ is contained in the union of supports of the corresponding operators $U_2^{\mb,j}$, and those supports do not overlap. Then the conclusion of Theorem $\ref{th_main_abstract}$ holds.
\end{cor}
\section{Particular cases and refinements of Theorem \ref{th_main_abstract}}
As the title suggests, in this section we will discuss several more concrete examples of the construction from Theorem \ref{th_main_abstract}. In each particular case, one can make some optimizations in choosing the unique continuation sets and obtain weaker requirements on the separation between different blocks.
\subsection{Example 1: the case $d=1$}
We will illustrate the above procedure, first, on the case where $f$ has a single flat segment and is $C_{\reg}$-regular outside of a neighborhood of it. Assume the following.
\begin{itemize}
	\item[(z1)] $f(x)=E$ for $x\in [a,a+L]=[b-L/2,b+L/2]\subset (-1/2,1/2)$.
	\item[(z2)] $L$ is not an integer multiple of $\omega$. Let $L=p\omega+z$, $p\in \N$, $z\in (0,1)$. Let also $\beta=\min\{z,1-z\}$.
	\item[(z3)] $f$ is $C_{\reg}$-regular outside of $[a-\beta,a+L+\beta]$.
	\item[(z4)] Let $M=\lceil L/2\omega \rceil+2$. Assume that the points
	$$
	  b-M\omega,b-(M-1)\omega,\ldots,b,b+\omega,\ldots,b+(M+1)\omega.
	$$
	are contained inside the interval $(-1/2,1/2)$.
\end{itemize}
Condition $(z4)$ is the most important one and has to be present in some form. It manifests the separation condition (2) from Theorem \ref{th_main_abstract}. Let also
$$
E_{\reg}:=\max\{|f(b-M\omega)|,|f(b+(M+1)\omega)|\}.
$$
\begin{thm}
\label{th_main_1d}
Fix $C_{\reg}$, $C_{\dio}$, $\tau_{\dio}$, $D_{\min}=1$, $\beta>0$, $M\in \N$. Let also $\omega\in \mathrm{DC}_{C_{\dio},\tau_{\dio}}\cap(0,1/2)$. Let $f$ satisfy $(f1)$ and $(z1)$ -- $(z4)$ with the $M$, $\beta$ fixed previously. There exists
$$
\ep_0=\ep_0(C_{\reg}, C_{\dio}, \tau_{\dio},\beta,M,E_{\reg})>0
$$
such that, for $0<\ep<\ep_0$, the operator \eqref{eq_h_def} is unitarily equivalent to an operator with convergent perturbation series. As a consequence, $H(x)$ satisfies Anderson localization for all $x\in \R$.
\end{thm}
\begin{proof}
Let $H_{M}(x)$ be the restriction of $H(x)$ into the interval $[-M,\ldots,M+1]$, a $(2M+2)\times (2M+2)$-matrix. The asymmetry is caused by the fact that we will consider $H(x)$ for $x\in [b,b+\omega]$.

For $x\in [b,b+\omega]$, consider the matrix $H_M(x)$ and modify it as follows: the entries that couple $-M$ with $-M+1$ will be linearly interpolated from $\ep$ at $x=b$ to $0$ at $x=b+\omega$, and the entries coupling $M$ and $M+1$, respectively, in the opposite way: $0$ at $x=b$ and $\ep$ at $x=0$. Denote the resulting operator $H_M'(x)$. Denote also by $U_M(x)$ the diagonalising operator obtained in Theorem \ref{th_diagonalisation1} for $H_{M}'(x)$. Proposition \ref{prop_jacobi_separation} combined with Proposition \ref{prop_epsilon2} guarantees that the family $U_M(x)$ satisfies the assumptions of Theorem \ref{th_diagonalisation1} (note that we cannot apply \ref{prop_jacobi_separation} directly, since, after rescaling, the range of energies will be of the size $\ep^{-1}$; however, one can use \ref{prop_epsilon2} to treat large eigenvalues corresponding to near the edges of the interval).

Denote by $U^{\mb}_0(x)$ the matrix $U_M(x)$ extended as identity into $\ell^2(\Z)$. At the moment, it has only been defined for $x\in [b,b+\omega]$. However, the block structure of $U^{\mb}_0(x)$ implies (for a single value of $x=b$ for which the equality makes sense):
\bee
\label{eq_covariant_1d}
U^{\mb}_0(x+\omega)=T^{-1}U^{\mb}_0(x)T.
\ene
As in the previous section, we will use \eqref{eq_covariant_1d} as the definition of $U_{0}^{\mb}(x)$ for $x\in \R$. We now perform the above procedure with each ``copy'' of the singular interval appearing in $H(x)$. That is, define
$$
 U_1^{\mb}(x)=\prod_{m\in\Z}U^{\mb}_0(x+m).
$$
Note that (z4) guarantees that the supports of the factors in the above product are disjoint. We do not need to introduce $U_2^{\mb}$, since there are no other lattice directions.

Now, let
$$
 H_1(x)=U^{\mb}_1(x)^{-1}H(x)U^{\mb}_1(x).
$$
At $x=a$, the diagonal matrix elements of $H(a)$ at $[-M+3,M-3]$ correspond to the values of $f$ calculated on the flat segment. They are exactly equal to the eigenvalues of $H'(a)$. In the notation of the previous section, we have $S=[-M+3,M-3]$. 

As $x$ goes from $a$ to $a+\omega$, the position of the flat segment on $\Z$ shifts by $1$: the entry located at $-M+3$ becomes regular, and the one located at $M-2$ becomes singular, and the flat segment now covers the interval $[-M+2,M-2]$. The entries at $-M,-M+1,-M+2,M-1,M,M+1$ are always regular and are small perturbations of the entries of $H(x)$, and therefore inherit their monotonicity properties. The remaining eigenvalues correspond to the set $[-M+3,M+2]$, which will play the role of $S\cup (S+\be_1))$ from the previous section:
$$
S\cup (S+\be_1)=[-M+3,M-2].
$$
For each $x\in [b,b+\omega]$, consider the points 
\bee
\label{eq_lattice_points}
x-M\omega,x-(M-1)\omega,\ldots,x,x+\omega,\ldots,x+(M+1)\omega.
\ene
Assumption (z2) guarantees that both endpoints of $[a,a+L]$ cannot be close to the points \eqref{eq_lattice_points} at the same time. More precisely, either $a$ or $a+L$ is at least $\beta/2$ away from the points \eqref{eq_lattice_points}. As a consequence, one can pick two adjacent points $p,p+1\in [-M,M+1]$ such that $x+p\omega$ and $x+(p+1)\omega$ are two closest points to $[a,a+L]$ on one side, but still not too close (both are $\beta/2$-away). For example, if $x=b$, then the points $p=-M+1$, $p+1=-M+2$ would work. As $x$ increases, we may have to switch to a pair of points on the other side. The set $\{p,p+1\}$ will play the role of $B$, the unique continuation set from the previous section. As discussed previously, this set is $x$-dependent; however, all cases can be treated similarly, and we will assume $B=\{-M+1,-M+2\}$ for simplicity.

 Let $\psi$ be an eigenfunction corresponding to one of the lattice points on $S\cup (S+\be_1)=[-M+3,M-2]$, normalized in, say, $\ell^2([-M,M+1])$, with eigenvalue $E_j(x)$, $j\in [-M+3,M-2]$ (here, the eigenvalues are identified with lattice points through Theorem \ref{th_diagonalisation1}). By Proposition \ref{prop_clustering}, $\psi$ is supported on $[-M+3,M-2]$ up to $O(\varepsilon)$. Moreover, $\psi(-M+2)=O(\varepsilon)$, $\psi(\-M+1)=O(\varepsilon^2)$. The values $f(x+m\omega)$, $m\in S$, are all equal to $E$, except maybe for $f(x+(M-2)\omega)$.
 
 We will now discuss the unique continuation procedure from $B=\{-M+1,-M+2\}$ into $[-M+3,M-2]$. Proposition \ref{prop_jacobi_separation} combined with Proposition \ref{prop_clustering} imply that
 $$
 E_{j+1}(x)-E_j(x)\ge c_1(\omega,M)\ep.
 $$
 On the other hand, by considering the eigenvalues as perturbations of the diagonal elements, we also have
$$
|E_j(x)-E|\le 3\ep,\,j\in [-M+3,M-3];\quad |E_{M-2}(x)-f(x+(M-2)\omega)|\le 3\ep.
$$
In both cases, the right hand side is actually $2\ep+c_2(\omega,M)\ep^2$, and we assumed $\ep$ to be small enough for simplicity. Assume first that $\psi$ is a singular eigenfunction corresponding to $E_j$ with $j\in [-M+3,M-3]$. By iterating the eigenvalue equation, it is easy to see that 
\bee
\label{eq_unique_M2}
|\psi(-M+2)|\ge c(M,\omega,E_{\reg})\ep.
\ene
Indeed, if $|\psi(-M+2)|$ is too small, then since $|\psi(-M+1)|=O(\ep^2)$, iterating the eigenvalue equation will contradict the $\ell^2$-normalization condition). Apply Proposition \ref{prop_hell} and conclude:
$$
E_j'(x)\ge c(M,\omega,E_{\reg})\ep^2,\quad j\in [-M+3,M-3].
$$
For the remaining eigenvalue $E_{M-3}$, we need to consider two cases. If, say, $|E-f(x+(M-2)\omega)|<10\ep$, then the previous argument also applies to the eigenfunction corresponding to $E_{M-3}$. If $|E-f(x+(M-2)\omega)|\ge 10\ep$, then the point $M-2$ can be considered as an individual cluster, and the corresponding eigenfunction must be non-trivially supported at that point (one can check that $|\psi(M-2)|\ge \frac12 \|\psi\|_{\ell^2})$. Then, one can apply Proposition \ref{prop_hell} directly at that point and obtain an even stronger conclusion: 
$$
E_{M-2}'(x)\ge c(M,\omega,E_{\reg})>0.
$$
In other words, such point can be considered as a regular point for practical purposes.

The above argument implies (conv2) from Proposition \ref{prop_convergence} with $\mu=2$ for our operator. We will illustrate (conv3)  on the example with $M=5$. The result of applying Proposition \ref{prop_clustering} to the eigenvectors of $H'_M(x)$ can be described as follows:
$$
U_M(x)\approx\begin{pmatrix}
1&\ep&\ep^2&\ep^3&\ep^3&\ep^3&\ep^3&\ep^3&\ep^3&\ep^9&\ep^{10}&\ep^{11}\\
\ep&1& \ep& \ep^2 & \ep^2& \ep^2& \ep^2& \ep^2& \ep^2& \ep^8&\ep^{9}&\ep^{10}\\
\ep^{2}&\ep& 1& \ep & \ep& \ep& \ep& \ep& \ep& \ep^7& \ep^8&\ep^{9}\\
\ep^{3}&\ep^2&\ep& 1& 1& 1& 1& 1& 1& \ep^6&\ep^7&\ep^{8}\\
\ep^{4}&\ep^3&\ep^2& 1& 1& 1& 1& 1& 1& \ep^5&\ep^6&\ep^{7}\\
\ep^{5}&\ep^4&\ep^3& 1& 1& 1& 1& 1& 1& \ep^4&\ep^5&\ep^{6}\\
\ep^{6}&\ep^5&\ep^4& 1& 1& 1& 1& 1& 1& \ep^3&\ep^4&\ep^{5}\\
\ep^{7}&\ep^6&\ep^5& 1& 1& 1& 1& 1& 1& \ep^2&\ep^3&\ep^{4}\\
\ep^{8}&\ep^7&\ep^6& 1& 1& 1& 1& 1& 1& \ep&\ep^2&\ep^{3}\\
\ep^{9}&\ep^8& \ep^7& \ep & \ep& \ep& \ep& \ep& \ep& 1& \ep&\ep^2\\
\ep^{10}&\ep^9 & \ep^8& \ep^2 & \ep^2& \ep^2& \ep^2& \ep^2& \ep^2& \ep& 1&\ep\\
\ep^{11}&\ep^{10}&\ep^9&\ep^3&\ep^3&\ep^3&\ep^3&\ep^3&\ep^3&\ep^2&\ep&1
\end{pmatrix}
$$
An entry $\ep^k$ means that the corresponding component of the eigenvector is bounded above by $c\ep^k$, where $c$ does not depend on $\ep$. The above matrix illustrates the worst possible case for $x\in [a,a+\omega]$: for example, at least one entry in edge of the central block is, actually, $\ep$ instead of $1$, and some of the remaining entries in that column are actually better by one order of $\ep$. Remark \ref{prop_epsilon_loss} implies that the derivatives in $x$ of each eigenvector are at most $\ep$ worse than the actual components. The central $6\times 6$ block corresponds to the interval $[-2,3]$, which covers the locations of singular diagonal entries corresponding to the flat piece.

Since $U_M(x)$ diagonalises a block of $H'(x)$ rather than $H(x)$, the conjugation of $H(x)$ by $U^{\mb}_0(x)$ will result in some extra off-diagonal entries. To calculate them, let us subtract $H'(x)$ from $H(x)$ and calculate the result of the conjugation of the remaining part:
$$
\begin{pmatrix}
1&0&0\\
0&U_M(x)^T&0\\
0&0&1
\end{pmatrix}\cdot
\begin{pmatrix}
O(1)&\ep&0&0&0&0&0\\
\ep&0&\ep&0&0&0&0\\
0&\ep&0&0_{1\times 8}&0&0&0\\
0&0&0_{8\times 1}&0_{8\times 8}&0_{8\times 1}&0&0\\
0&0&0&0_{1\times 8}&0&\ep&0\\
0&0&0&0&\ep&0&\ep\\
0&0&0&0&0&\ep&O(1)\\
\end{pmatrix}\cdot
\begin{pmatrix}
1&0&0\\
0&U_M(x)&0\\
0&0&1
\end{pmatrix}.
$$
Each element of the matrix product is made out of products of three entries, one from each matrix. Even without multiplying three $(14\times 14)$-matrices, one can see the following: in order to produce a non-zero element that starts from an entry inside the $[-2,3]$ block, one needs to use a matrix entry from the right factor that connects point entry from $[-2,3]$ to a point outside of the interval $[-3,4]$, thus gaining at least $\ep^2$. If we gained exactly $\ep^2$, then one will gain at least one $\ep$ from the central factor. Thus, any jump that starts inside of the $[-2,3]$ block will cost at least $\ep^3$, which is sufficient to satisfy (conv3) and apply Proposition \ref{prop_convergence} (in view of (conv4) and Remark \ref{prop_epsilon_loss}).
\end{proof}
\begin{rem}
One can check that, under the assumptions of Theorem \ref{th_main_1d}, the derivative of the integrated density of states of the family $H(x)$ has a spike of height $\asymp\ep^{-2}$ and width $\asymp\ep^2$ around energy $E$. One can produce larger spikes by combining multiple flat intervals at different energies, see Examples 5 and 6 below.
\end{rem}
\subsection{A higher-dimensional version}
The approach from the previous subsection can be extended to several higher-dimensional examples. 

\vskip 2mm

{\noindent \bf Example 2:} Let $d\ge 2$. Suppose that $f$ and $\omega_1$ satisfy all assumptions from the previous subsection (with $\omega=\omega_1$). Additionally, assume the following
\begin{itemize}
	\item[(z5)]For any $x\in [a,a+L]$ and any $\bn\in \Z^d$, $\bn\notin \Z \be_{1}$, $|\bn|_1\le 6$, we have $f$ being $C_{\reg}$-regular at $x+\bn\cdot\omega$.
\end{itemize}
\begin{thm}
\label{th_main_zd}
Suppose that the assumptions of Theorem $\ref{th_main_1d}$ are satisfied with $\omega$ replaced by $\omega_1$. Assume also that $\omega\in\mathrm{DC}_{C_{\dio},\tau_{\dio}}$ and $(z5)$. Then there exists
$$
\ep_0=\ep_0(d,C_{\reg}, C_{\dio}, \tau_{\dio}, \beta,M,E_{\reg})>0
$$
such that the operator \eqref{eq_h_def} satisfies Anderson localization on $\Z^d$.
\end{thm}
\begin{proof}
We only outline the changes in the proof. Instead of the interval $[-M,M+1]$, consider the box
$$
[-3,3]^{d-1}\times[-M,M+1]
$$
and construct $U_M(x)$ by applying Theorem \ref{th_diagonalisation1} to the above box. One will need to interpolate the Laplacian terms that couple $[-3,3]^{d-1}\times\{-M\}$ and $[-3,3]^{d-1}\times\{-M+1\}$ between $\ep$ and $0$, and the terms coupling $[-3,3]^{d-1}\times\{M\}$ and $[-3,3]^{d-1}\times\{M+1\}$ between $0$ and $\ep$, as $x$ goes from $b$ to $b+\omega_1$. The eigenvalues on the ``central'' box $\{0\}^{d-1}\times [-M,M+1]$ are $c\ep$-separated due to Proposition \ref{prop_jacobi_separation}. Eigenvalues corresponding to the remaining entries are separated from each other by a constant independent of $\ep$ for $\ep\ll 1$ due to Proposition \ref{prop_clustering}. Therefore, expansion of the box in the  directions perpendicular to $\be_1$ will shift the eigenvalues of the block $\{0\}^{d-1}\times [-M,M+1]$ by at most $O(\ep^2)$ due to Proposition \ref{prop_epsilon2}, and therefore will preserve $O(\ep)$-separation property.

The rest follows from the proof of Theorem \ref{th_main_1d}. One needs to specify the unique continuation subset. For $x=a$ it will be
$$
B=\left(\{0\}^{d-1}\times\{-M+1,-M+2\}\right)\cup \left(\{\bn'\in \Z^{d-1}\colon |\bn'|_1=1\}\times [-M+1,M-3]\right).
$$
In other words, we are allowed to make up to two steps to the left from the singular subset, and one step in one other direction. The set $B$ can be imagined as a ``claw'' surrounding the singular set $S$.

Proposition \ref{prop_clustering} implies that any singular eigenfunction restricted to $\{-M+1\}\times \{\bn'\in \Z^{d-1}\colon |\bn'|_1=1\}$ will be $O(\ep^2)$. The rest of the set $B$ is one step away from $S$, and therefore any eigenfunction will have, ultimately, at least $c\ep$ of mass somewhere on $B$. One may have to treat the right-most point of $S$ separately, the argument does not change from Theorem \ref{th_main_1d}, in view of Proposition \ref{prop_jacobi_separation} and Proposition \ref{prop_epsilon2}.

Once $U_M(x)$ is constructed for $x\in [b,b+\omega_1]$, we repeat the steps in Section 5 in order to construct the operator $U_2^{\mb}$. Condition (z3) will guarantee that the supports of different copies of $U_0^{\mb}(x)$ will not overlap.
\end{proof}
\subsection{Several additional cases}
In this subsection, we outline several additional cases which, ultimately, are covered by Theorem \ref{th_main_abstract}.

\vskip 2mm

{\noindent \bf Example 3:} One can consider multiple intervals of the type described in Theorem \ref{th_main_1d}, as long as they are separated enough so that the supports of the operator $U_1$ are disjoint. Similarly, one can combine multiple ``non-interacting'' cases of Theorem \ref{th_main_zd}.

\vskip 2mm

{\noindent \bf Example 4:} More interestingly, if two intervals from Theorem \ref{th_main_1d} are not sufficiently disjoint, one can consider them as one singular set and apply the argument from the general theorem. One can also do it in the setting of Theorem \ref{th_main_zd}, with a strengthened version of (z5).

\vskip 2mm

{\noindent \bf Example 5:} The following example can be called ``a chain of intervals''. Let $d=1$ and Let $I_1,\ldots,I_N\subset (-1/2,1/2)$ be a collection of disjoint intervals: 
$$
	f(x)=E_j\,\,\,\text{for}\,\,x\in I_j;\quad E_j< E_k\,\,\,\text{for}\,\,j<k,
$$
and for $j=1,\ldots,N-1$. Assume also that, for sufficiently large $r=r(N)$ ($r=2N$ should be enough), we have 
$$
-1/2<I_1-r\omega<I_{N}+r\omega<1/2.
$$
Then the assumptions of Theorem \ref{th_main_abstract} are satisfied, and we have Anderson localization for small $\ep$.
\begin{rem}
Suppose that, on top of all these assumptions, that $I_{j+1}=I_j+\omega$; and we also assume for simplicity that $N=2k-1$. Then it is easy to see that $f_2'(x)\asymp \ep^{2k}$ for $x\in I_k$. As a result, the derivative of the integrated density of states has a spike of height $\asymp \ep^{-2k}$ and width $\asymp\ep^{2k}$ around $E_k$.
\end{rem}

{\noindent \bf Example 6:} The following example is a combination of Examples 2 and 5. It can be called ``a tree of intervals''. Let $I_1,\ldots,I_N$ be a collection of disjoint intervals: 
	$$
	f(x)=E_j\,\,\,\text{for}\,\,x\in I_j;\quad E_j\neq E_k\,\,\,\text{for}\,\,j\neq k.
	$$
	Suppose that for some $x\in I_1$ we have, say, $x_2:=x+\omega_2\in I_2$. Additional shifts by $\omega_2,\ldots,\omega_d$ may produce points, say, $x_3=x_2+\omega_3\in I_3,\ldots$. Assume that, for any $x\in I_1$, we are guaranteed to escape all intervals after, say, $q$ steps in the directions perpendicular to $\be_1$, and would only be able to come back to $I_1$ after additional $r$ steps. In other words, for $x\in I_1$, $f$ is regular at $x+\omega'\cdot\bn'$, for all $\bn'\in \Z^{d-1}$, $q\le |\bn'|\le q+r$. One can refine this in each particular case, if the number of steps needed to escape the system of intervals depends on the direction. In all cases, we would require a ``layer'' of thickness $r$ that consists of regular points and surrounds the collection of singular points described above. The thickness of the layer will depend on the number of steps required to escape the singular intervals. For $r$ sufficiently large, the assumptions of Theorem \ref{th_main_abstract} will be satisfied. In particular, to establish the separation condition (5) from Theorem \ref{th_main_abstract}, note that each interval produces a block with $c\ep$-separated eigenvalues due to Proposition \ref{prop_jacobi_separation}. Different intervals are located at different energies, and therefore do not affect each other (each lies in its own cluster). The coupling between these intervals satisfies the assumptions of Proposition \ref{prop_epsilon2}, and therefore also does not affect the separation condition. After the conjugation, we will have, for the new diagonal function $f_2(x)$,
$$
f_2'(x)\asymp \ep^{\mu_j},\quad x\in I_j,
$$
where $\mu_j$ depends on the number of steps in the directions perpendicular to $\be_1$ required to escape the interval (by more careful analysis, one can show that it is equal to twice the number of steps). As a consequence, the derivative of the integrated density of states of $H$  will have spikes of order $\ep^{-\mu_j}$ around the energies close to $E_j$.
\end{document}